# THE CALCULATION OF EXPECTATIONS FOR CLASSES OF DIFFUSION PROCESSES BY LIE SYMMETRY METHODS

BY MARK CRADDOCK AND KELLY A. LENNOX

*University of Technology, Sydney*

This paper uses Lie symmetry methods to calculate certain expectations for a large class of Itô diffusions. We show that if the problem has sufficient symmetry, then the problem of computing functionals of the form $E_x(e^{-\lambda X_t - \int_0^t g(X_s)\,ds})$ can be reduced to evaluating a single integral of known functions. Given a drift $f$ we determine the functions $g$ for which the corresponding functional can be calculated by symmetry. Conversely, given $g$, we can determine precisely those drifts $f$ for which the transition density and the functional may be computed by symmetry. Many examples are presented to illustrate the method.

**1. Introduction.** Suppose that a stochastic process $X = \{X_t : t \geq 0\}$ is a solution of the Itô SDE

$$(1.1) \qquad dX_t = a(X_t, t)\,dt + b(X_t, t)\,dW_t, \qquad X_0 = x.$$

Here $W = \{W_t : t \geq 0\}$ is a standard Wiener process. We would like to compute the expectations

$$(1.2) \qquad \xi_{\lambda, g}(X_t) = E(e^{-\lambda X_t - \int_0^t g(X_s)\,ds} | X_0 = x).$$

These types of functionals have numerous applications. For example, if $\lambda = 0$, $g(x) = x$, then the expectation gives the price of a zero coupon bond in financial mathematics. If $g(X_t) = \mu h(X_t)$ and $h(X_t)$ does not depend on $\mu$, then the expectation is the two dimensional Laplace transform of the joint density of $(X_t, \int_0^t h(X_s)\,ds)$.

Most cases must be handled numerically, but recently considerable attention has been devoted to the problem of finding classes of diffusions for which the transition density and expectations such as (1.2) may be determined exactly. The papers [4, 5, 7, 8, 14] provide some details.

By the Feynman–Kac formula, if the solution of the Cauchy problem

$$
\begin{aligned}
u_t = \sigma x^\gamma u_{xx} + f(x)u_x - g(x)u, \qquad x > 0, \\
u(x,0) = \varphi(x), \qquad x \in \Omega = [0,\infty),
\end{aligned}
\tag{1.3}
$$

is unique, then it is given by

$$
u(x,t) = E_x(e^{-\int_0^t g(X_s)\,ds}\varphi(X_t)),
\tag{1.4}
$$

where $X_t$ satisfies the SDE $dX_t = f(X_t)\,dt + \sqrt{2\sigma X_t^\gamma}\,dW_t$. We will use Lie symmetry methods to solve (1.3) by finding an explicit fundamental solution of the PDE as a function of the drift. We have a number of results which do this. The first is Theorem 3.1 below.

The Cauchy problem (1.3) does not usually have a unique solution unless we impose certain extra conditions (see, e.g., Part 4 of [18]). In fact, our methods actually yield more than one fundamental solution for the operator in (1.3). See Section 3. The selection of the correct fundamental solution for the purposes of computing the expectation must therefore be treated with care.

The plan is as follows. In Section 2, for a given function $g$, we show that the drift $f$ must satisfy one of a family of Riccati equations if the PDE is to have nontrivial Lie symmetries. We also show how to explicitly solve the Riccati equations. In Section 3, we derive classes of generalized Laplace transforms of our fundamental solutions. Our methods require a time independent solution of the PDE and we show how to obtain these. Sections 4 and 5 contain examples and general formulae for the fundamental solutions. We also use group invariant solution methods to compute the precise form of the fundamental solutions in cases where the transforms are difficult to invert.

In the final section, given a drift $f$, we determine those functions $g$ for which the Cauchy problem (1.3) can be solved by symmetry and the expectation $\xi_{\lambda,g}(X_t)$ can computed.

Each method for determining exact solutions has its own unique strengths and weaknesses. Specific techniques work well for particular class of problems, but may fail when we move outside that class. For example, methods which use eigenfunction expansions and stochastic analysis have been exploited successfully for many years. The literature is enormous, but see [14] for some recent work. However useful descriptions of the necessary eigenfunctions are only available for limited classes of PDEs and typically the eigenfunction expansion of a solution cannot be explicitly summed. Eigenfunction methods also become more difficult when we move to higher dimensions.

Lie symmetry methods have natural advantages. They are easy to use, cover a broad range of cases and one can develop very simple tests to decide whether or not a PDE is amenable to symmetry analysis. The major limitation is that many important PDEs do not have any useful symmetries.

However, for the rich class of PDEs amenable to Lie group analysis, the technique is extremely effective. We obtain explicit closed form expressions for the fundamental solutions in terms of the drift and do not need any changes of variables or measure. A major advantage is that the method can be extended to higher dimensions. An understanding of Lie techniques in one dimension is essential if we wish to apply them in higher dimensional problems.

**2. Calculating the Lie symmetries.** A *symmetry* of a differential equation is a transformation which maps solutions to solutions. In the 1880s Lie developed a technique for systematically determining all groups of *point symmetries* for systems of differential equations.[1] Olver's book [15], is a fine modern account of Lie's theory of symmetry groups as are the texts by Bluman and Kumei [2], Hydon [12] and Ovsiannikov [16]. For a survey of some recent work on Lie symmetries, see [6].

In [7] Craddock and Lennox proved that if a PDE of the form $u_t = \sigma x^\gamma u_{xx} + f(x)u_x - \mu x^n u$ has nontrivial symmetries, then there is always a symmetry arising from an action of $SL(2,\mathbb{R})$ that is a classical integral transform of a fundamental solution. We generalize the results of Craddock and Lennox and apply them to problems in stochastic calculus.

We will work with PDEs of the form (1.3). We observe that every PDE of the form (1.3) possesses time translation symmetries as well as scaling symmetries in $u$. That is, $u(x, t + \varepsilon)$ is a solution if $u(x, t)$ is and so is $cu(x, t)$ for any constant $c$. We term these trivial symmetries. We require that the PDE have nontrivial symmetries.

PROPOSITION 2.1. *If $\gamma \neq 2$, the PDE (1.3) has a nontrivial Lie symmetry group if and only if $f$ is a solution of one of the following families of drift equations.*

$$Lf = Ax^{2-\gamma} + B, \tag{2.1}$$

$$Lf = Ax^{2-\gamma} + Bx^{1-\gamma/2} - \frac{3}{8(2-\gamma)}\sigma^2, \tag{2.2}$$

*where*

$$Lf = \sigma x^\gamma \left(\frac{x^{1-\gamma}f(x)}{2\sigma(2-\gamma)}\right)'' + f(x)\left(\frac{x^{1-\gamma}f(x)}{2\sigma(2-\gamma)}\right)' + g(x) + \frac{xg'(x)}{2-\gamma}. \tag{2.3}$$

---

[1] There also exist so called *generalized symmetries*, but we will not consider them.

*If $\gamma = 2$ then the PDE has a nontrivial Lie group of symmetries if and only if*

(2.4) $$Uf = A,$$

(2.5) $$Uf = A \ln x + B.$$

With $v(x) = \frac{f(x) \ln x}{x}$ we have

$$Uf = \frac{x^2}{4}v''(x) + \frac{f(x)}{4\sigma}v'(x) - \frac{f(x)}{4x} + \frac{xg'(x)\ln x}{2} + g(x).$$

PROOF. The proof of this result is textbook, so we just summarize it for the $\gamma \neq 2$ case. By Lie's method, we look for infinitesimal symmetries of the form

(2.6) $$\mathbf{v} = \xi(x,t,u)\partial_x + \tau(x,t,u)\partial_t + \phi(x,t,u)\partial_u,$$

where we employ the convention $\frac{\partial}{\partial x} = \partial_x$ etc. We seek to determine conditions on $\xi, \tau$ and $\phi$ which guarantee that $\mathbf{v}$ generates a symmetry of the PDE. Standard arguments show that $\tau$ can only depend on $t$ and $\xi$ can only depend on $x$ and $t$. Further, $\phi$ must be linear in $u$. Lie's Theorem says that $\mathbf{v}$ generates a local group of symmetries if and only if $\text{pr}^2 \mathbf{v}[u_t - (\sigma x^\gamma u_{xx} + f(x)u_x - g(x)u)] = 0$ whenever $u$ is a solution of (1.3). Here $\text{pr}^2 \mathbf{v}$ is the second prolongation of $\mathbf{v}$. (See Chapter 2 of [15] for the explicit prolongation formula for a vector field.) If $\gamma \neq 2$, this leads to the conditions

(2.7) $$\xi(x,t) = \frac{x}{2-\gamma}\tau_t + x^{\gamma/2}\rho, \qquad \phi(x,t,u) = \alpha(x,t)u + \beta(x,t),$$

(2.8) $$\alpha(x,t) = \frac{-x^{2-\gamma}}{2\sigma(2-\gamma)^2}\tau_{tt} - \frac{x^{1-\gamma/2}}{\sigma(2-\gamma)}\rho_t + \left(\frac{\gamma}{4}x^{\gamma/2-1} - \frac{x^{-\gamma/2}}{2\sigma}f(x)\right)\rho$$
$$- \frac{x^{1-\gamma}}{2\sigma(2-\gamma)}f(x)\tau_t + \eta.$$

Here $\tau, \rho, \eta$ are functions of $t$ only and satisfy

(2.9) $$\frac{-x^{2-\gamma}}{2\sigma(2-\gamma)^2}\tau_{ttt} - \frac{x^{1-\gamma/2}}{\sigma(2-\gamma)}\rho_{tt} + \eta_t = \frac{-(1-\gamma)}{2(2-\gamma)}\tau_{tt} - Lf\tau_t + Kf\rho,$$

with $Kf = \sigma x^\gamma(\frac{\gamma}{4}x^{\gamma/2-1} - \frac{x^{-\gamma/2}}{2\sigma}f(x))'' + f(x)(\frac{\gamma}{4}x^{\gamma/2-1} - \frac{x^{-\gamma/2}}{2\sigma}f(x))$. The function $\beta(x,t)$ is an arbitrary solution of (1.3). The proof is completed when we compare the corresponding powers of $x$ and recognize that $\tau$ will only be a nontrivial function of $t$ if $f$ is a solution of one of the given equations. The $\gamma = 2$ case is similar. $\square$

2.1. *Solving the drift equations.* If we set $h(x) = x^{1-\gamma} f(x)$ then the ODE $Lf = Ax^{2-\gamma} + B$ is equivalent to the Riccati equation

(2.10) $\quad \sigma x h' - \sigma h + \frac{1}{2}h^2 + 2\sigma x^{2-\gamma} g(x) = \sigma A x^{4-2\gamma} + 2\sigma B x^{2-\gamma} + C.$

Similarly for the second of this pair of drift equations. These Riccati equations can be linearized by the change of variables $h = 2\sigma x y'/y$.

We point out that it is natural to consider two different cases with equation (2.10), namely $A = 0, A \neq 0$, because the Lie symmetries of the PDE take a quite different form if $A = 0$.

Finding solutions of the drift equations in the $\gamma = 2$ case is simplified by the substitution $H(\ln x) = \frac{f(x) \ln x}{x} - \sigma \ln x$. Under which, for example, $Uf = A$ becomes

$$\sigma \xi H' - \sigma H + \frac{1}{2}H^2 + 2\sigma \xi^2 g(\xi) = \left(2A + \frac{\sigma}{2}\right)\xi^2 + B.$$

Thus we have precisely the same form of Riccati equation.

2.2. *The $\gamma = 1$ case.* If $\gamma = 1$ the drift equations are

(2.11) $\quad \sigma x f' - \sigma f + \frac{1}{2}f^2 + 2\sigma x g(x) = Ax + B,$

(2.12) $\quad \sigma x f' - \sigma f + \frac{1}{2}f^2 + 2\sigma x g(x) = \frac{1}{2}Ax^2 + Bx + C,$

(2.13) $\quad \sigma x f' - \sigma f + \frac{1}{2}f^2 + 2\sigma x g(x) = \frac{1}{2}Ax^2 + \frac{2}{3}Bx^{3/2} + Cx - \frac{3}{8}\sigma^2.$

Here $A, B, C$ are arbitrary constants. The factors of $\frac{1}{2}$ and $\frac{2}{3}$ multiplying $A$ and $B$ are a notational convenience. We will concentrate in this paper on the first two Riccati equations, which we treat separately. The third case is harder and will be discussed elsewhere.

**3. Applications of the symmetries.** Symmetries map solutions to solutions. If the symmetry group is sufficiently rich, then one can construct extremely complex solutions of a given PDE from trivial solutions. We are able to obtain an integral transform of the fundamental solution by applying a symmetry to a suitable stationary solution. We begin with the $A = 0$ case.

THEOREM 3.1. *Suppose that $h(x) = x^{1-\gamma} f(x)$ is a solution of the Riccati equation*

(3.1) $\quad \sigma x h' - \sigma h + \frac{1}{2}h^2 + 2\sigma x^{2-\gamma} g(x) = 2\sigma A x^{2-\gamma} + B.$

*Then the PDE*

(3.2) $\quad u_t = \sigma x^\gamma u_{xx} + f(x) u_x - g(x) u, \qquad x \geq 0,$

*has a symmetry of the form*

$$\overline{U}_\varepsilon(x,t) = \frac{1}{(1+4\varepsilon t)^{(1-\gamma)/(2-\gamma)}} \exp\left\{\frac{-4\varepsilon(x^{2-\gamma} + A\sigma(2-\gamma)^2 t^2)}{\sigma(2-\gamma)^2(1+4\varepsilon t)}\right\}$$

$$\times \exp\left\{\frac{1}{2\sigma}\left(F\left(\frac{x}{(1+4\varepsilon t)^{2/(2-\gamma)}}\right) - F(x)\right)\right\}$$

$$\times u\left(\frac{x}{(1+4\varepsilon t)^{2/(2-\gamma)}}, \frac{t}{1+4\varepsilon t}\right),$$

*where $F'(x) = f(x)/x^\gamma$ and $u$ is a solution of the PDE. That is, $U_\varepsilon$ is a solution of (3.2) whenever $u$ is. If $u(x,t) = u_0(x)$ with $u_0$ an analytic, stationary solution then there is a fundamental solution $p(t,x,y)$ of (3.2) such that*

$$(3.3) \qquad \int_0^\infty e^{-\lambda y^{2-\gamma}} u_0(y) p(t,x,y)\, dy = U_\lambda(x,t).$$

*Here $U_\lambda(x,t) = \overline{U}_{(1/4)\sigma(2-\gamma)^2 \lambda}$.*

Equation (3.3) is a generalized Laplace transform. We may recover the fundamental solution by inverting the transform and dividing by $u_0(y)$. Before proceeding to the proof we need a technical result.

PROPOSITION 3.2.  *The solution $U_\lambda(x,t)$ in Theorem 3.1 is the Laplace transform of a distribution.*

PROOF. This follows from the observation that $U_\lambda(x,t)$ can be written as a product of $\lambda^\nu$ for some value $\nu$ and an analytic function $G(1/\lambda)$. Any function which is analytic in $1/\lambda$ is a Laplace transform. Further $\lambda^\nu$ is the Laplace transform of a distribution. The product of two Laplace transforms is a Laplace transform. Hence $U_\lambda(x,t)$ is a Laplace transform. See the table in [19], page 348 for the inverse Laplace transform of $\lambda^\nu$ for different values of $\nu$ and Chapter 10 for general results on when a distribution can be represented as an inverse Laplace transform. □

Now we prove Theorem 3.1.

PROOF OF THEOREM 3.1. Lie's method shows that (3.2) has an infinitesimal symmetry of the form

$$(3.4) \quad \mathbf{v} = \frac{8xt}{2-\gamma}\partial_x + 4t^2\partial_t - \left(\frac{4x^{2-\gamma}}{\sigma(2-\gamma)^2} + \frac{4x^{1-\gamma}tf(x)}{\sigma(2-\gamma)} + \beta t + 4At^2\right)u\partial_u,$$

where $\beta = \frac{4(1-\gamma)}{2-\gamma}$. Exponentiating this symmetry and applying it to a solution $u(x,t)$ yields $\overline{U}_\varepsilon$. The idea is that a solution $U_\lambda(x,t)$ of the PDE should be obtained as

$$U_\lambda(x,t) = \int_0^\infty U_\lambda(y,0)p(t,x,y)\,dy.$$

The symmetry solution has the property that $U_\lambda(x,0) = e^{-\lambda x^{2-\gamma}}u_0(x)$. Thus we should have

$$\int_0^\infty e^{-\lambda y^{2-\gamma}}u_0(y)p(t,x,y)\,dy = U_\lambda(x,t).$$

But this is the generalized Laplace transform of $u_0 p$. We thus need to show that $U_\lambda$ is a generalized Laplace transform of some distribution $u_0 p$. Since

$$\int_0^\infty e^{-\lambda y^{2-\gamma}}u_0(y)p(t,x,y)\,dy$$
$$= \int_0^\infty e^{-\lambda z}u_0(z^{1/(2-\gamma)})p(t,x,z^{1/(2-\gamma)})\frac{z^{1/(2-\gamma)-1}\,dz}{2-\gamma},$$

we must show that $U_\lambda$ is the Laplace transform of some distribution $u_0 p$. This follows from Proposition 3.2.

To see that $p$ is a fundamental solution of the PDE, observe that if we integrate a test function $\varphi(\lambda)$ with sufficiently rapid decay against $U_\lambda$ then the function $u(x,t) = \int_0^\infty U_\lambda(x,t)\varphi(\lambda)\,d\lambda$ is a solution of (3.2). We also have

$$u(x,0) = \int_0^\infty U_\lambda(x,0)\varphi(\lambda)\,d\lambda = \int_0^\infty u_0(x)e^{-\lambda x^{2-\gamma}}\varphi(\lambda)\,d\lambda = u_0(x)\Phi(x),$$

where $\Phi$ is the generalized Laplace transform of $\varphi$. Next observe that by Fubini's Theorem

$$\int_0^\infty u_0(y)\Phi(y)p(t,x,y)\,dy = \int_0^\infty \int_0^\infty u_0(y)\varphi(\lambda)p(t,x,y)e^{-\lambda y^{2-\gamma}}\,d\lambda\,dy$$
$$= \int_0^\infty \int_0^\infty u_0(y)\varphi(\lambda)p(t,x,y)e^{-\lambda y^{2-\gamma}}\,dy\,d\lambda$$
$$= \int_0^\infty \varphi(\lambda)U_\lambda(x,t)\,dx = u(x,t).$$

We know that $u(x,0) = u_0(x)\Phi(x)$. Thus integrating initial data $u_0\Phi$ against $p$ solves the Cauchy problem for (3.2), with this initial data. Hence $p$ is a fundamental solution. $\square$

COROLLARY 3.3. *If $g = 0$ in (1.3), then there is a fundamental solution $p(t,x,y)$ with the property that*

$$\int_0^\infty e^{-\lambda y^{2-\gamma}}p(t,x,y)\,dy$$

$$(3.5) \qquad = \frac{1}{(1+4\varepsilon t)^{(1-\gamma)/(2-\gamma)}} \exp\left\{\frac{-4\varepsilon(x^{2-\gamma} + A\sigma(2-\gamma)^2 t^2)}{\sigma(2-\gamma)^2(1+4\varepsilon t)}\right\}$$

$$\times \exp\left\{\frac{1}{2\sigma}\left(F\left(\frac{x}{(1+4\varepsilon t)^{2/(2-\gamma)}}\right) - F(x)\right)\right\},$$

and $\int_0^\infty p(t,x,y)\,dy = 1$. Here $\varepsilon = \frac{1}{4}\sigma(2-\gamma)^2\lambda$.

PROOF. Since $g(x) = 0$ we may take $u_0(x) = 1$ in Theorem 3.1. Observe that $U_0(x,t) = 1$. Thus $\int_0^\infty p(t,x,y)\,dy = 1$. $\square$

For the $\gamma = 2$ case we have a similar result.

THEOREM 3.4. *Suppose that*

$$(3.6) \qquad \frac{x^2}{4}\left(\frac{f(x)\ln x}{x}\right)'' + \frac{f(x)}{4\sigma}\left(\left(\frac{f(x)\ln x}{x}\right)' - \frac{\sigma}{x}\right) + \frac{x\ln x g'(x)}{2} + g(x) = A.$$

*Let $u_0(x)$ be a stationary solution of*

$$(3.7) \qquad u_t = \sigma x^2 u_{xx} + f(x) u_x - g(x) u$$

*analytic near zero. Then there is a fundamental solution of (3.7) such that*

$$(3.8) \qquad \int_0^\infty e^{-(\varepsilon/\sigma)(\ln y)^2} p(t,x,y) u_0(y)\,dy = U_\varepsilon(x,t),$$

*where*

$$(3.9) \qquad U_\varepsilon(x,t) = \frac{1}{\sqrt{1+4\varepsilon t}}\exp\left\{-\frac{\varepsilon((\ln x)^2 + 2\sigma t\ln x + (4A+\sigma)\sigma)t^2)}{\sigma(1+4\varepsilon t)}\right\}$$

$$\times \exp\left\{\frac{1}{2\sigma}\left(F\left(\frac{\ln x}{1+4\varepsilon t}\right) - F(x)\right)\right\} u_0(x^{1/(1+4\varepsilon t)})$$

*and $F'(x) = e^{-x} f(x)$.*

PROOF. The proof is similar to the previous result. If the drift $f$ satisfies (3.6), then the PDE has an infinitesimal symmetry of the form

$$(3.10) \qquad \mathbf{v} = 4xt\ln x\frac{\partial}{\partial x} + 4t^2\frac{\partial}{\partial t}$$

$$-\left(\frac{(\ln x)^2}{\sigma} + \ln x^2\left(\frac{f(x)}{\sigma x} - 1\right)t + 2t + (4A+\sigma)t^2\right)u\frac{\partial}{\partial u}.$$

Exponentiating $\mathbf{v}$ and applying it to $u_0$ produces $U_\varepsilon(x,t)$. This has the initial value $U_\varepsilon(x,0) = u_0(x)e^{-(\varepsilon/\sigma)(\ln x)^2}$. The proof follows the same lines as previously. $\square$

### 3.1. Finding stationary solutions.

For $\gamma \neq 2$, to obtain a stationary solution $u_0(x)$ given $h$ we need to solve $\sigma x^\gamma u_{xx} + f(x)u_x - g(x)u = 0$. We divide by $x^{\gamma-1}$ to rewrite this as $\sigma x u_{xx} + h(x)u_x - x^{1-\gamma}g(x)u = 0$. If $u(x) = \tilde{u}(x)e^{-(1/(2\sigma))\int h(x)/x\,dx}$, then $2\sigma^2 x^2 \tilde{u}''(x) - (2\sigma A x^{2-\gamma} + B)\tilde{u}(x) = 0$. Finally, the substitution $z = x^{2-\gamma}$, $\tilde{u}(x) = w(x^{2-\gamma})$ reduces this ODE to

$$2\sigma^2(2-\gamma)^2 z^2 w''(z) + 2\sigma^2(2-\gamma)(1-\gamma)zw'(z) - (2\sigma Az + B)w(z) = 0.$$

This equation has Bessel function solutions. Thus, given $f$ and $g$ we can always determine the integral transform (3.3) and hence express the fundamental solution, at least up to a Laplace inversion integral.

### 3.2. Finding different fundamental solutions.

Different choices of $u_0$ will in general lead to *different* fundamental solutions. In fact we may sometimes extract two fundamental solutions from a single choice of $u_0$. We illustrate the first situation with an example involving a squared Bessel process of dimension 3. The PDE $u_t = 2xu_{xx} + 3u_x$ has two stationary solutions, $u_0(x) = 1$ and $u_0(x) = 1/\sqrt{x}$. The choice $u_0(x) = 1$ leads to

$$(3.11) \qquad \int_0^\infty e^{-\lambda y} p(t,x,y)\,dy = \frac{1}{(1+2\lambda t)^{3/2}} \exp\left(-\frac{\lambda x}{1+2\lambda t}\right),$$

which is easily inverted giving

$$(3.12) \qquad p(t,x,y) = \frac{1}{\sqrt{2\pi t x}} e^{-(x+y)/(2t)} \sinh\left(\frac{\sqrt{xy}}{t}\right).$$

This is the transition density for the squared Bessel process of dimension 3. If we use the second stationary solution we get

$$(3.13) \qquad \int_0^\infty \frac{1}{\sqrt{y}} e^{-\lambda y} p_2(t,x,y)\,dy = \frac{1}{\sqrt{x(1+2\lambda t)}} \exp\left(-\frac{\lambda x}{1+2\lambda t}\right).$$

From which we obtain a second fundamental solution

$$(3.14) \qquad p_2(t,x,y) = \frac{1}{\sqrt{2\pi t x}} e^{-(x+y)/(2t)} \cosh\left(\frac{\sqrt{xy}}{t}\right).$$

This is not a transition density. Whereas $\int_0^\infty p(t,x,y)\,dy = 1$ we have

$$(3.15) \qquad \int_0^\infty p_2(t,x,y)\,dy = \sqrt{\frac{2t}{\pi x}} e^{-x/(2t)} + \mathrm{erf}\left(\sqrt{\frac{x}{2t}}\right).$$

The function (3.15) is a solution of $u_t = 2xu_{xx} + 3u_x$ which satisfies $\lim_{t \to 0} u(x,t) = 1$ for all $x > 0$, but it is not continuous at the origin. All solutions generated by $p_2(t,x,y)$ have this feature. They solve the Cauchy problem for the PDE on any domain $[\varepsilon, \infty)$, $\varepsilon > 0$, but are not continuous at the origin.

Next consider the equation

(3.16) $$u_t = xu_{xx} + \left(3 - \frac{4b}{b+ax^2}\right)u_x, \qquad x > 0.$$

An application of Corollary 3.3 yields the existence of a fundamental solution $p(t,x,y)$ satisfying

(3.17) $$\int_0^\infty e^{-\lambda y} p(t,x,y)\,dy = \frac{ax^2 + b(1+\lambda t)^4}{(b+ax^2)(1+\lambda t)^3} \exp\left\{-\frac{\lambda x}{1+\lambda t}\right\}.$$

We take $a, b$ positive. To invert this Laplace transform requires the introduction of distributions. (A complete discussion of the inversion of the Laplace transforms which require distribution theory is found in [9].) Inverting the Laplace transform produces the fundamental solution

(3.18) $$\begin{aligned}p(t,x,y) &= \frac{x}{yt}\frac{b+ay^2}{b+ax^2}e^{-(x+y)/t}I_2\left(\frac{2\sqrt{xy}}{t}\right)\\ &\quad + \frac{b(x+t)e^{-(x+y)/t}}{t(b+ax^2)}\delta(y) + \frac{bte^{-x/t}}{b+ax^2}\delta'(y).\end{aligned}$$

Here $\delta'(y)$ is derivative of the Dirac delta function, $I_\nu$ is the modified Bessel function of the first kind. It is easy to show that $\int_0^\infty p(t,x,y)\,dy = 1$. If $q(t,x,y) = p(t,x,y) - \frac{bte^{-x/t}}{b+ax^2}\delta'(y)$, then $q$ is also a fundamental solution and $\int_0^\infty q(t,x,y)\,dy = 1$ as well, since $\delta'$ acts by $\int_0^\infty \phi(y)\delta'(y)\,dy = -\phi'(0)$. Indeed we may also take $r(t,x,y) = \frac{x}{yt}\frac{b+ay^2}{b+ax^2}e^{-(x+y)/t}I_2\left(\frac{2\sqrt{xy}}{t}\right)$ and this is a fundamental solution with

(3.19) $$\int_0^\infty r(t,x,y)\,dy = 1 - e^{-x/t}\frac{b(t+x)}{t(b+ax^2)} = l(x,t).$$

Now $l(x,t) \to 1$ as $t \to 0$ for all $x \geq 0$. Thus (3.16) with the data $u(x,0) = 1$ does not have a unique solution. If we wish to interpret $p(t,x,y)$ and $q(t,x,y)$ as probabilities we are faced with a situation in which we have two densities arising from the same generator that define the same null set and probabilities, but because of the terms involving generalized function they do not define the same expectations.

The drift function $f$ has the property that $-1 \leq f(x) < 3$. This means that if we consider in a formal sense the SDE $dX_t = f(X_t)\,dt + \sqrt{2X_t}\,dW_t$ we are in a situation where $X_t$ may become negative because of the drift and hence the diffusion $X_t$ may then become complex valued. It is not obvious how to interpret this. The more natural SDE $dX_t = f(X_t)\,dt + \sqrt{2|X_t|}\,dW_t$ leads to the PDE $u_t = 2|x|u_{xx} + f(x)u_x$, but our theory does not cover a PDE of this form, though one can also analyse PDEs of this type by symmetry. A full discussion of the problem of trying to associate a diffusion with a PDE of the form (1.3) covered by our methods is beyond the scope of this paper. But it is clear that our techniques allow us to construct very general classes of PDEs which have very complex behavior.

For our purposes we will restrict attention to problems where the interpretation of the fundamental solutions as densities is clear. That is, we will look at generators which arise from well behaved SDEs with positive drifts and/or unique strong, real valued solutions. We will include examples of fundamental solutions with delta function terms, but our examples will be more straightforward than the above. A full investigation of the phenomenon we have described here will be considered elsewhere.

This leaves the question of which choice of $u_0$ we take in order to use the Feynman–Kac formula to compute expectations. The examples we give will show how this choice is made.

EXAMPLE 3.1. Consider the Bessel process $X = \{X_t : t \geq 0\}$ where

$$dX_t = \frac{a}{X_t}\,dt + dW_t, \qquad X_0 = x. \tag{3.20}$$

Bessel processes are well studied and the transition density for this process is known. Here $\sigma = 1/2, \gamma = 0$ and $f(x) = \frac{a}{x}$. We will require $a > \frac{1}{2}$ for simplicity. So $h(x) = a$ and we easily see that $A = 0$ in the Riccati equation. Set $g(x) = \frac{\mu}{4x^2}$. We want a fundamental solution of

$$u_t = \frac{1}{2}u_{xx} + \frac{a}{x}u_x - \frac{\mu}{4x^2}u, \qquad \mu > 0.$$

If we are to calculate $E_x(e^{-(\mu/4)\int_0^t \frac{ds}{X_s^2}}\phi(X_t))$ then we require a fundamental solution that reduces to the known transition density of the process at $\mu = 0$.

A stationary solution is given by solving $2x^2 u'' + 4ax u' - \mu u = 0$. We take the solution $u_0(x) = x^d$ where $d = \frac{1}{2} - a + \sqrt{\frac{\mu}{2} + (a - \frac{1}{2})^2}$. Observe that as $\mu \to 0$, $u_0(x) \to 1$. Setting $\varepsilon = \lambda\sigma = \frac{\lambda}{2}$ in (3.3) gives

$$\int_0^\infty e^{-\lambda y^2} y^d p(t,x,y)\,dy = \frac{x^d}{(1+2\lambda t)^\nu}\exp\left(-\frac{\lambda x^2}{1+2\lambda t}\right), \tag{3.21}$$

with $\nu = d + a + \frac{1}{2}$. This is a Laplace transform if we let $y^2 = z$. Inverting the transform yields

$$p(t,x,y) = \frac{y}{t}\left(\frac{y}{x}\right)^{a-1/2}\exp\left(-\frac{x^2+y^2}{2t}\right)I_{\nu-1}\left(\frac{xy}{t}\right). \tag{3.22}$$

If $\mu = 0$, this reduces to the well known transition density of a Bessel process. We can now write $E_x(e^{-\lambda X_t - (\mu/4)\int_0^t \frac{ds}{X_s^2}}) = \int_0^\infty e^{-\lambda y}p(t,x,y)\,dy.$

This gives the Laplace transform of the joint density of $(X_t, \int_0^t \frac{ds}{X_s^2})$. We cannot evaluate this integral. However

$$E_x(e^{-\lambda X_t^2 - (\mu/4)\int_0^t \frac{ds}{X_s^2}}) = \int_0^\infty e^{-\lambda y^2} p(t,x,y)\,dy$$

$$= e^{-x^2/(2t)}\left(\frac{x^2}{2t}\right)^{(2\nu - 2a - 1)/4}$$

$$\times \frac{\Gamma(\alpha)\,_1F_1(\alpha, \nu, x^2/(2t + 4t^2\lambda))}{\Gamma(\nu)(1 + 2t\lambda)^\alpha},$$

where $\alpha = \frac{1+2a+2\nu}{4}$ and $_1F_1(a,b,z)$ is Kummer's confluent hypergeometric function.

There is another solution $x^{d_2}$ with $d_2 = \frac{1}{2} - a - \sqrt{\frac{\mu}{2} + (a - \frac{1}{2})^2}$. We observe that as $\mu \to 0$ $x^{d_2} \to x^\nu$ where $\nu = \frac{1}{2} - a - |a - \frac{1}{2}|$ which is nonconstant. So the fundamental solution coming from this stationary solution does not reduce to the necessary density. This situation is typical of these problems.

NOTE. There are many Laplace transforms which can be computed by other methods. However our methods often yield different representations of these quantities than can be found in much of the literature. For a Bessel process of index $\xi$ there is a well known alternative representation

$$(3.23) \quad E_x(e^{-\lambda X_t^2 - (\mu^2/2)\int_0^t \frac{ds}{X_s^2}})$$
$$= \frac{x^{\gamma - \xi}}{\Gamma(p)} \int_0^\infty \frac{v^{p-1} \exp(-x^2(v + \lambda)/(1 + 2(v + \lambda)t))}{(1 + 2(v + \lambda)t)^{1+\gamma}}\,dv,$$

where $p = \frac{1}{2}(\gamma - \xi)$ and $\gamma^2 = \xi^2 + \mu^2$. This is proved using the absolute continuity law $P_x^{(\xi)}|_{\mathcal{F}_t} = (\frac{X_t}{x})^\xi e^{-(\mu^2/2)\int_0^t \frac{ds}{X_s^2}} P_x^{(0)}|_{\mathcal{F}_t}$. However we have never seen the integral (3.23) explicitly evaluated. Our representation allows for the explicit evaluation of the integral from tables. Moreover to obtain the joint density, we only have to invert a one dimensional Laplace transform rather than a two dimensional Laplace transform, since we already know the inverse Laplace transform in the $\lambda$ variable.

EXAMPLE 3.2. Now we turn to the Bessel process with drift. The SDE is

$$dX_t = \left(\frac{a + 1/2}{X_t} + \frac{bI_{a+1}(bX_t)}{I_a(bX_t)}\right)dt + dW_t, \quad a > -1.$$

These processes arise from the radial part of $d$-dimensional Brownian motion with drift and were studied by Pitman and Yor in [17]. Here $\sigma = \frac{1}{2}$ and $h = xf$ satisfies $\frac{1}{2}xh' - \frac{1}{2}h + \frac{1}{2}h^2 = \frac{b^2}{2}x^2 + \frac{1}{8}(4a^2 - 1)$. So $A = \frac{1}{2}b^2$. Using the stationary solution $u_0(x) = 1$ we obtain

$$(3.24) \quad \int_0^\infty e^{-\lambda y^2} p(t,x,y)\,dy = \frac{I_a(bx/(1+2\lambda t))}{(1+2\lambda t)I_a(bx)} \exp\left(-\frac{\lambda(x^2 + b^2 t^2)}{1+2\lambda t}\right).$$

We convert this to a Laplace transform by setting $y^2 = z$. Inverting the Laplace transform gives the transition density

$$(3.25) \quad p(t,x,y) = \frac{y}{t}\frac{I_a(by)}{I_a(bx)}\exp\left(-\frac{x^2+y^2}{2t} - \frac{b^2}{2}t\right)I_a\left(\frac{xy}{t}\right).$$

Since $\frac{I_a(by)}{I_a(bx)} \to (\frac{y}{x})^a$ as $b \to 0$, this reduces to the transition density of a Bessel process as $b \to 0$. Now consider the problem of computing $E_x(e^{-\lambda X_t - \mu \int_0^t \frac{ds}{X_s^2}})$. We require a stationary solution of the PDE

$$(3.26) \quad u_t = \frac{1}{2}u_{xx} + \left(\frac{a+1/2}{x} + \frac{bI_{a+1}(bx)}{I_a(bx)}\right)u_x - \frac{\mu}{x^2}u.$$

We look for a solution of the form $u_0(x) = e^{-F(x)}v(x)$ where $F' = f$. We see that

$$4x^2 v''(x) + 4(1 - 4a^2 - 4b^2x^2 - 8\mu)v(x) = 0.$$

Since $F(x) = \ln(\sqrt{x}I_a(bx))$, we obtain the stationary solution $u_0^\mu(x) = \frac{I_{\sqrt{a^2+2\mu}}(bx)}{I_a(bx)}$. Note that as $\mu \to 0$, $u_0^\mu \to 1$. Theorem 3.1 tells us that there is a fundamental solution $q(t,x,y)$ of (3.26) such that

$$(3.27) \quad \int_0^\infty e^{-\lambda y^2} u_0^\mu(y) q(t,x,y)\,dy \\ = \frac{e^{-\lambda(b^2 t^2 + x^2)/(1+2t\lambda)} I_{\sqrt{a^2+2\mu}}(bx/(1+2t\lambda))}{(1+2t\lambda)I_a(bx)}.$$

Inversion gives the fundamental solution

$$(3.28) \quad q(t,x,y) = \frac{y}{t}\frac{I_a(by)}{I_a(bx)}\exp\left(-\frac{x^2+y^2}{2t} - \frac{b^2}{2}t\right)I_{\sqrt{a^2+2\mu}}\left(\frac{xy}{t}\right).$$

Clearly this reduces to the transition density at $\mu = 0$. Obviously $E_x(e^{-\lambda X_t - \mu \int_0^t \frac{ds}{X_s^2}}) = \int_0^\infty e^{-\lambda y} q(t,x,y)\,dy$, but we have not been able to evaluate this integral. However to find the joint density, we only need to invert the Laplace transform $q(t,x,y)$ in $\mu$. This can be done numerically.

The procedure we have is now clear. We begin with an SDE which possesses a unique solution and compute a transition density for it. This can be done by our methods. We choose the fundamental solution of (1.3) which reduces to this density when $g \to 0$. From this we may compute the desired expectations.

**4. Laplace transforms of joint densities.** As is clear from Example 3.1, it is possible to use the techniques we have developed here to compute Laplace transforms of joint densities very easily. We will focus attention on the $\gamma = 1$ case and do an example of the $\gamma = 0$ case below. It should be clear that we may prove similar results for any $\gamma$. If $\gamma = 1, g(x) = \mu x^n$, with $n = \pm 1$ and $n = \frac{1}{2}$, then we can find drift functions which do not depend on $\mu$. It is possible to easily compute Laplace transforms of the joint densities of $(X_t, \int_0^t X_s)$, $(X_t, \int_0^t \frac{ds}{X_s})$ and indeed we can also find the Laplace transform of the joint density of $(X_t, \int_0^t X_s \, ds, \int_0^t \frac{ds}{X_s})$. The four dimensional Laplace transforms of the law of $(X_t, \int_0^t X_s \, ds, \int_0^t \frac{ds}{X_s}, \int_0^t \sqrt{X_s} \, ds)$ can be computed for certain classes of diffusions, but this is more complicated and will be treated elsewhere.

EXAMPLE 4.1. Let $X = \{X_t : t \geq 0\}$ be a squared Bessel process, where
$$dX_t = n \, dt + 2\sqrt{X_t} \, dW_t, \qquad X_0 = x.$$
To compute $u(x,t) = E_x(e^{-\mu \int_0^t \frac{ds}{X_s}} h(X_t))$, we require a fundamental solution for the PDE

(4.1) $$u_t = 2x u_{xx} + n u_x - \frac{\mu}{x} u.$$

This is the case when $\sigma = 2, \gamma = 1, g(x) = \mu/x$ in Theorem 2.1. The stationary solutions are $u_0(x) = x^d$, $d_\pm = \frac{1}{4}(2 - n \pm \sqrt{(n-2)^2 + 8\mu})$. We take $d = d_+$, as $x^{d_-}$ diverges near zero and does not tend to the constant solution as $\mu \to 0$. Applying the symmetry we have

(4.2) $$\int_0^\infty y^d p(t,x,y) e^{-\lambda y} \, dy = \frac{x^d}{(1+2\lambda t)^{2d+n/2}} \exp\left\{\frac{-\lambda x}{1+2\lambda t}\right\}.$$

This Laplace transform may immediately be inverted to give

(4.3) $$p(t,x,y) = \frac{1}{2t}\left(\frac{x}{y}\right)^{(1-n/2)/2} I_{2d+n/2-1}\left(\frac{\sqrt{xy}}{t}\right) \exp\left\{-\frac{(x+y)}{2t}\right\}.$$

From this we obtain the two dimensional Laplace transform of the joint density of $(X_t, \int_0^t \frac{1}{X_s} \, ds)$.

$$\int_0^\infty e^{-\lambda y} p(t,x,y) \, dy = e^{-x/2t} \left(\frac{x}{2t}\right)^d \frac{\Gamma(\alpha) {}_1F_1(\alpha, \beta, x/(2t+4t^2\lambda))}{\Gamma(\beta)(1+2\lambda t)^\alpha},$$

with $\alpha = d + \frac{n}{2}, \beta = 2d + \frac{n}{2}$.

We have not been able to invert this Laplace transform. It is interesting to compare this result when $\lambda = 0$ with the famous Hartman–Watson law for a squared Bessel process which gives the *conditional* expectation $E_x(e^{-(\mu^2/2)\int_0^t \frac{ds}{X_s}}|X_t = y) = \frac{I_{\sqrt{\mu^2+\nu^2}}(\sqrt{xy}/t)}{I_\nu(\sqrt{xy}/t)}$, where $\nu = 1 - \frac{1}{2}n$.

EXAMPLE 4.2. Now we turn to the SDE

$$(4.4) \qquad dX_t = \frac{aX_t}{1 + aX_t/2} dt + \sqrt{2X_t}\, dW_t, \qquad X_0 = x, a > 0.$$

For the PDE $u_t = xu_{xx} + \frac{ax}{1+ax/2}u_x - \frac{\mu}{x}u$ we have two stationary solutions $u_1(x) = \frac{x^{(1+\sqrt{1+4\mu})/2}}{2+ax}$ and $u_2(x) = \frac{x^{(1-\sqrt{1+4\mu})/2}}{2+ax}$, neither of which is equal to 1 at $\mu = 0$. However, the linear combination $2u_2 + au_1$ does satisfy the desired condition, so we use the stationary solution $u_0(x) = x^{(1-\sqrt{1+4\mu})/2}\left(\frac{2+ax\sqrt{1+4\mu}}{2+ax}\right)$. Letting $d_\pm = \frac{1}{2}(1 \pm \sqrt{1+4\mu})$ we obtain

$$(4.5) \quad U_\lambda(x,t) = \frac{1}{2+ax}\left(\frac{ax^{d_+}}{(1+\lambda t)^{2d_+}} + \frac{2x^{d_-}}{(1+\lambda t)^{2d_-}}\right)\exp\left(-\frac{\lambda x}{1+\lambda t}\right).$$

Inverting the Laplace transform yields the fundamental solution

$$(4.6) \qquad \begin{aligned} p(t,x,y) &= \frac{1}{t}\frac{2+ay}{2+ax}\sqrt{\frac{x}{y}}\frac{ay^{d_+}}{2y^{d_+}+ay^{d_+}}I_{2d_+-1}\left(\frac{2\sqrt{xy}}{t}\right)e^{-(x+y)/t} \\ &\quad + \frac{2+ay}{t^{2d_-}(2+ax)}\frac{2x^{d_-}}{2y^{d_+}+ay^{d_+}} \\ &\quad \times e^{-(x+y)/t}\int_0^y I_0\left(\frac{2\sqrt{x\xi}}{t}\right)\nu_\mu(\xi - y)\,d\xi. \end{aligned}$$

Here $\nu_\mu(\xi) = \mathcal{L}^{-1}[\lambda^{\sqrt{1+4\mu}}]$. This inverse Laplace transform is a right sided distribution. See [19] for a table of these inverse Laplace transforms. We must be careful when taking the limit, however a rather laborious calculation shows that as $\mu \to 0$ this reduces to the fundamental solution

$$(4.7) \qquad p(t,x,y) = \frac{e^{-(x+y)/t}}{(2+ax)t}\left[\sqrt{\frac{x}{y}}(2+ay)I_1\left(\frac{2\sqrt{xy}}{t}\right) + t\delta(y)\right],$$

which is the transition density for the process $X_t$. From (4.6) we may calculate $E_x(e^{-\lambda X_t - \mu \int_0^t \frac{ds}{X_s}})$. The integral can be evaluated exactly, but we leave this to the reader.

We may readily establish more general results. We present an example shortly. First however, we need a lemma.

LEMMA 4.1. *The PDE $u_t = \sigma x u_{xx} + f(x) u_x - \frac{\mu}{x} u$, with $f$ a solution of the Riccati equation $\sigma x f' - \sigma f + \frac{1}{2} f^2 = Ax + B$ has a stationary solution $u_0^\mu(x)$ which has the property that $u_0^0(x) = 1$.*

PROOF. This is a direct calculation. Setting $f = 2\sigma x \frac{y'}{y}$ reduces the Riccati equation to $2\sigma x^2 y'' - (Ax + B) y = 0$ which has general solution $y = \sqrt{x}(c_1 I_\alpha(\sqrt{2Ax}) + c_2 K_\alpha(\sqrt{2Ax}))$, where $\alpha = \frac{1}{\sigma}\sqrt{2B + \sigma^2}$ and $K_\alpha$ is the usual modified Bessel function. Now a stationary solution of the PDE is given by solving $\sigma x u'' - f(x) u' - \frac{\mu}{x} u = 0$. The substitution $u_0(x) = e^{-F(x)/(2\sigma)} v(x)$ leads to the ODE $2\sigma^2 x^2 v'' - (Ax + B + 2\mu\sigma) v = 0$, where $F'(x) = f(x)/x$. Thus there is a stationary solution of the PDE

$$(4.8) \qquad u_0^\mu(x) = \frac{c_1 I_\nu(\sqrt{2Ax}) + c_2 K_\nu(\sqrt{2Ax})}{c_1 I_\alpha(\sqrt{2Ax}) + c_2 K_\alpha(\sqrt{2Ax})},$$

where $\nu = \frac{1}{\sigma}\sqrt{2B + \sigma^2 + 4\mu\sigma}$. Clearly $u_0^0(x) = 1$. □

THEOREM 4.2. *Suppose that $f$ is a solution of the Riccati equation $\sigma x f' - \sigma f + \frac{1}{2} f^2 = Ax + B$ with $A > 0$. Then the PDE*

$$u_t = \sigma x u_{xx} + f(x) u_x - \frac{\mu}{x} u, \qquad 2B + \sigma^2 + 4\mu\sigma > 0$$

*has fundamental solution*

$$(4.9) \qquad \begin{aligned} p(t,x,y) &= \frac{\sqrt{x}}{\sigma t} e^{-F(x)/2\sigma - (x+y)/(\sigma t) - (A/(2\sigma)) t} \frac{1}{u_0^\mu(y)} \\ &\quad \times \left( c_1 I_\nu\left(\frac{2\sqrt{xy}}{\sigma t}\right) I_\nu\left(\frac{\sqrt{Ay}}{\sigma}\right) \right. \\ &\qquad \left. + c_2 I_{-\nu}\left(\frac{2\sqrt{xy}}{\sigma t}\right) I_{-\nu}\left(\frac{\sqrt{Ay}}{\sigma}\right) \right), \end{aligned}$$

*where $F'(x) = f(x)/x$, $u_0^\mu(x)$ is as above and $\nu = \frac{1}{\sigma}\sqrt{2B + \sigma^2 + 4\mu\sigma}$ satisfies $|\nu| < 1$. From which we may calculate $E(e^{-\lambda X_t - \mu \int_0^t \frac{ds}{X_s}})$ by*

$$(4.10) \qquad E_x(e^{-\lambda X_t - \mu \int_0^t \frac{ds}{X_s}}) = \int_0^\infty e^{-\lambda y} p(t,x,y) \, dy.$$

PROOF. Let $f(x) = 2\sigma x y'/y$ where $y(x) = \sqrt{x}(c_1 I_\alpha(\frac{\sqrt{2Ay}}{\sigma}) + c_2 I_{-\alpha}(\frac{\sqrt{2Ay}}{\sigma}))$. (Since $\alpha$ is not an integer $K_\alpha(z) = I_\alpha(z)$.) To obtain a stationary solution of the PDE we use as in the previous lemma the substitution $u_0(x) = e^{-F(x)/(2\sigma)} v(x)$ in the equation $\sigma x u'' + f(x) u' - \mu/x u = 0$. Then $v$ must satisfy $2\sigma x^2 v''(x) - (Ax + B + 2\sigma\mu) v(x) = 0$, which has solutions $\sqrt{x} I_{\pm\nu}(\frac{\sqrt{2Ax}}{\sigma})$.

Set
$$u_0^\mu(x) = \sqrt{x}e^{-F(x)/(2\sigma)}\left(c_1 I_\nu\left(\frac{\sqrt{2Ax}}{\sigma}\right) + c_2 I_{-\nu}\left(\frac{\sqrt{2Ax}}{\sigma}\right)\right).$$

We then have $u_0^0(y) = 1$. Now

(4.11)
$$\int_0^\infty e^{-\lambda y} u_0^\mu(y) p(t,x,y)\, dy$$
$$= \frac{\sqrt{x}e^{-F(x)/(2\sigma)-\lambda(x+At^2/2)/(1+\lambda\sigma t)}}{1+\lambda\sigma t}$$
$$\times \left[c_1 I_\nu\left(\frac{\sqrt{2Ax}}{\sigma(1+\lambda\sigma t)}\right) + c_2 I_{-\nu}\left(\frac{\sqrt{2Ax}}{\sigma(1+\lambda\sigma t)}\right)\right].$$

We use the fact that
$$\mathcal{L}^{-1}\left[\frac{1}{\lambda}\exp\left(\frac{m^2+n^2}{\lambda}\right)I_d\left(\frac{2mn}{\lambda}\right)\right] = I_d(2m\sqrt{y})I_d(2n\sqrt{y}), \qquad d > -1.$$

Inverting the Laplace transform gives the fundamental solution. The expectation follows from the Feynman–Kac formula. $\square$

For $\mu$ outside the range of the theorem, the inversion of the Laplace transform will in general involve right sided distributions. Techniques for inverting these transforms are discussed in [9].

**5. The other Riccati equations.** Suppose that the drift $f$ is such that $h = x^{1-\gamma}f$ satisfies a Riccati equation of the form
$$\sigma x h' - \sigma h + \tfrac{1}{2}h^2 + 2\sigma x g(x) = Ax^{4-2\gamma} + Bx^{2-\gamma} + C, \qquad A \neq 0,$$
then we do not get Laplace type transforms of the fundamental solution, except in some special cases. Instead we get a more general integral transform of the fundamental solution, known as a Whittaker transform. A similar situation arises with the final Riccati equation, but we will consider that elsewhere as there are some technical issues which need to be carefully discussed.

The reason for the difference is that the Lie symmetries for these Riccati equations are different in form to the symmetries in the previous case. Consider the $\gamma = 1$ case. Craddock and Lennox proved in [7] that if $f$ is a solution of the second class of Riccati equations with $A \neq 0$ and $u_0(x)$ is a stationary solution of $u_t = \sigma x u_{xx} + f(x)u_x - \mu x^r u$, then so is

(5.1)
$$U_\varepsilon(x,t) = e^{-Bt/(2\sigma)}\exp\left\{\frac{-\sqrt{A}x\varepsilon}{2\sigma(e^{\sqrt{A}t} - \varepsilon)} + \frac{1}{2\sigma}\left(F\left(\frac{xe^{\sqrt{A}t}}{e^{\sqrt{A}t} - \varepsilon}\right) - F(x)\right)\right\}$$
$$\times (e^{\sqrt{A}t} - \varepsilon)^{B/(2\sigma\sqrt{A})} u_0\left(\frac{xe^{\sqrt{A}t}}{e^{\sqrt{A}t} - \varepsilon}\right).$$

In general this does not reduce to the necessary form for a Laplace transform when we set $t = 0$. Moreover, for certain stationary solutions, the action is trivial. However, there are special cases when it is a Laplace transform. We present an example.

EXAMPLE 5.1. Consider the PDE
$$u_t = xu_{xx} + \frac{ax}{1+ax/2}u_x - \mu x u.$$

The drift is a solution of the Riccati equation $\sigma x f' - \sigma f + \frac{1}{2}f^2 + 2\mu x^2 = \frac{1}{2}Ax^2 + Bx + C$, with $A = 4\mu, B = C = 0$. There are two stationary solutions, $u_0^{\pm}(x) = e^{\pm\sqrt{\mu}x}(ax+2)^{-1}$. The symmetry (5.1) leaves $u_0^+$ unchanged. However for the other solution we have

$$(5.2) \qquad U_\varepsilon(x,t) = \frac{e^{-\sqrt{\mu}x}}{2+ax}\exp\left(-\frac{2\sqrt{\mu}\varepsilon x}{e^{2\sqrt{\mu}t}-\varepsilon}\right).$$

By our general argument, with $\lambda = \frac{2\sqrt{\mu}\varepsilon}{1-\varepsilon}$ we can find a fundamental solution $p(t,x,y)$ such that

$$(5.3) \qquad \begin{aligned}\int_0^\infty e^{-\lambda y}\frac{e^{-\sqrt{\mu}y}}{2+ay}&p(t,x,y)\,dy \\ &= \frac{e^{-\sqrt{\mu}x}}{2+ax}\exp\left(\frac{-2\lambda\sqrt{\mu}x}{\lambda(e^{2\sqrt{\mu}t}-1)+2\sqrt{\mu}e^{2\sqrt{\mu}t}}\right).\end{aligned}$$

Inverting the Laplace transform yields the fundamental solution

$$\begin{aligned}p(t,x,y) &= \frac{2+ay}{2+ax}\frac{e^{\sqrt{\mu}y}}{e^{\sqrt{\mu}x}}\exp\left\{\frac{-2\sqrt{\mu}(x+ye^{2\sqrt{\mu}t})}{e^{2\sqrt{\mu}t}-1}\right\} \\ &\quad \times \left(\sqrt{\frac{\mu x}{y}}\frac{I_1(2\sqrt{\mu xy}/\sinh(\sqrt{\mu}t))}{\sinh(\sqrt{\mu}t)} + \delta(y)\right).\end{aligned}$$

Letting $\mu \to 0$ gives the transition density for the process $X$ with $dX_t = \frac{aX_t}{1+aX_t/2}dt + \sqrt{2X_t}\,dW_t$. Consequently we can compute the Laplace transform of the joint density of the processes $(X_t, \int_0^t X_s\,ds)$. Formula 6.621.3 of [11] and the Feynman–Kac formula give

$$\begin{aligned}E_x(e^{-\lambda X_t - \mu\int_0^t X_s\,ds}) &= \frac{e^{-\sqrt{\mu}x - 2\sqrt{\mu}x/(e^{2\sqrt{\mu}t}-1)}}{2+ax} \\ &\quad \times \left(2 + \frac{\sqrt{\mu}x(e^{(\beta^2/4\alpha)}(4/\beta + a\beta/(2\alpha^2)) - 4/\beta)}{\sinh(\sqrt{\mu}t)}\right).\end{aligned}$$

Here $\alpha = \lambda - \sqrt{\mu} + \frac{2\sqrt{\mu}e^{2\sqrt{\mu}t}}{e^{2\sqrt{\mu}t}-1}$ and $\beta = 2\sqrt{x\mu}\,\mathrm{cosech}(\sqrt{\mu}t)$.

Although we have recovered the necessary fundamental solution by inverting a Laplace transform, for these classes of Riccati equations we in general have to deal with a transform involving Whittaker functions.[2]

DEFINITION 5.1 (*The Whittaker transform*). The Whittaker transform of a suitable function $\phi$ is defined by

$$(5.4) \qquad \Phi(\lambda) = \int_0^\infty (\lambda y)^{-k-1/2} e^{-\lambda y/2} W_{k+1/2,\nu}(\lambda y) \phi(y) \, dy.$$

An inversion theorem for this transform is known. For a suitable constant $\rho$ we have, in the principal value sense

$$\phi(y) = \frac{1}{2\pi i} \frac{\Gamma(1+\nu-k)}{\Gamma(1+2\nu)} \int_{\rho-i\infty}^{\rho+i\infty} (\lambda y)^{-k-1/2} e^{\lambda y/2} M_{k-1/2,\nu}(\lambda y) \Phi(\lambda) \, d\lambda.$$

The functions $W_{k+1/2,\nu}$ and $M_{k-1/2,\nu}$ are the Whittaker functions given in 13.1 of [1]. For a discussion of the transform, see [3], page 110. Special cases of the Whittaker transform include the Laplace transform and Hankel transforms.

In [7] the following result is established for the $\gamma = 1, g(x) = \mu x^r$ case. The proof of the following result is identical to the case in [7].

THEOREM 5.2. *Let $f$ be a solution of the Riccati equation*

$$\sigma x f' - \sigma f + \tfrac{1}{2} f^2 + 2\sigma x g(x) = \tfrac{1}{2} A x^2 + B x + C.$$

*Let $\eta = \frac{B}{2\sigma\sqrt{A}} - \tfrac{1}{2}\beta$, $\beta = 1 + \sqrt{1 + \frac{2C}{\sigma^2}}$, $\nu = \tfrac{1}{2}\sqrt{1 + 2C/\sigma^2}$ and $k + \tfrac{1}{2} = -\frac{B}{2\sigma\sqrt{A}}$. Let $U_{\sqrt{A}/\lambda\sigma - 1}(x,t)$ be given by setting $\varepsilon \to \frac{\sqrt{A}}{\lambda\sigma} - 1$ in*

$$(5.5) \quad U_\varepsilon(x,t) = e^{-(Bt)/(2\sigma)} (e^{\sqrt{A}t} - \varepsilon)^{B/(2\sigma\sqrt{A}) - \beta/2} x^{\beta/2} e^{-F(x)/(2\sigma)}$$

$$\times \exp\left\{ \frac{-\sqrt{A}x(e^{\sqrt{A}t} + \varepsilon)}{2\sigma(e^{\sqrt{A}t} - \varepsilon)} \right\} \Psi\left(\alpha, \beta, \frac{\sqrt{A}xe^{\sqrt{A}t}}{\sigma(e^{\sqrt{A}t} - \varepsilon)}\right),$$

*where $\Psi(\alpha, \beta, z)$ is the Tricomi confluent hypergeometric function and $F'(x) = f(x)/x$. Suppose that $\lambda^{B/(\sigma\sqrt{A})} U_{\sqrt{A}/(\lambda\sigma)-1}(x,t)$ is the Whittaker transform of a function $\tilde{h}_\mu(t,x,y)$. Then this function $\tilde{h}_\mu$ is given by*

$$\tilde{h}_\mu(t,x,y) = \left(\frac{\sqrt{A}}{\sigma}\right)^\eta y^{k+1/2} e^{(\sqrt{A}y - F(y))/(2\sigma)} p_\mu(t,x,y),$$

*where $p_\mu$ is a fundamental solution of the PDE $u_t = \sigma x u_{xx} + f(x) u_x - g(x) u$.*

---

[2] The transform is due to Meijer and is often called a Meijer transform. However there are several other integral transforms which are also called Meijer transforms in the literature and this nomenclature can be somewhat confusing.

Similar results can be proved for the general case when $\gamma \neq 1$. The Whittaker transform is analytic, thus one can in principle invert a Whittaker transform as a series expansion, since the Whittaker transform of $y^n$ is of the form $C\lambda^{-n+1}$ where $C$ is a constant depending on $n$ and the parameters in the transform. Unfortunately, tables of Whittaker transforms tend to be very sparse, so in practice, these transforms are really only invertible in the special cases where they reduce to Laplace or Hankel transforms. Further, theorems guaranteeing that a given distribution is a Whittaker transform are available only for limited parameter ranges. At present therefore, these results are of mainly theoretical interest.

There are some other issues that arise when we try to calculate expectations of the form $E_x(e^{-\lambda X_t - \mu \int_0^t X_s \, ds})$. The stationary solutions of the PDE $u_t = \sigma x^\gamma u_{xx} + f(x)u_x - \mu x u$ do not necessarily have the property that a linear combination tends to the constant solution when $\mu \to 0$. See Example 5.1. Conversely both of the linearly independent stationary solutions may converge to a constant when we take the limit. We thus need a different criterion for deciding which fundamental solution to use.

An alternative approach to the construction of the necessary fundamental solutions is to use group invariant solutions. A group invariant solution is a solution which is left invariant under the action of a group transformation. Chapter three of Olver's book [15] contains a detailed discussion of group invariant solutions. See also the material in [6]. We will construct the general form of the fundamental solutions and determine the necessary one by examining the integrability properties.

The following result, which combines two theorems in Bluman and Kumei's text [2] is the key to our method.

THEOREM 5.3. *Consider the nth order boundary value problem $P(x, D^\alpha u) = 0$ subject to the conditions $B_j(x, u, u^{(n-1)}) = 0$ on the surface $\omega_j(x) = 0$. A vector field $\mathbf{v}$ is admitted by the BVP if:*

(1) $\mathrm{pr}^n \mathbf{v}[P(x, D^\alpha u)] = 0$ *when* $P(x, D^\alpha u) = 0$.
(2) $\mathbf{v}(\omega_j(x)) = 0$ *when* $\omega_j(x) = 0$.
(3) $\mathrm{pr}^{n-1} \mathbf{v}[B_j(x, u, u^{(n-1)})] = 0$ *when* $B_j(x, u, u^{(n-1)}) = 0$ *on the surface* $\omega_j(x) = 0$.

*Suppose that a BVP admits a vector field $\mathbf{v}$. Then the solution of the BVP is a group invariant solution with respect to the symmetries generated by $\mathbf{v}$.*

We prove several results as an illustration of the technique. We concentrate on the case $A > 0$, leaving the case $A < 0$ to the interested reader.

THEOREM 5.4. *Suppose that $f$ is a solution of the Riccati equation $\sigma x f' - \sigma f + \frac{1}{2}f^2 + 2\sigma\mu x^2 = \frac{1}{2}Ax^2 + Bx + C, A > 0$. Then there is a fundamental solution of the PDE*

$$u_t = \sigma x u_{xx} + f(x)u_x - \mu x u, \qquad x \geq 0,$$

*of the form*

$$p(t,x,y)$$
(5.6)
$$= \frac{\sqrt{Axy}e^{-(F(x)-F(y))/(2\sigma)}}{2\sigma \sinh(\sqrt{A}t/2)} \exp\left(-\frac{Bt}{2\sigma} - \frac{\sqrt{A}(x+y)}{2\sigma \tanh(\sqrt{A}t/2)}\right)$$
$$\times \left(C_1(y)I_\nu\left(\frac{\sqrt{Ayx}}{\sigma \sinh(\sqrt{A}t/2)}\right) + C_2(y)I_{-\nu}\left(\frac{\sqrt{Ayx}}{\sigma \sinh(\sqrt{A}t/2)}\right)\right),$$

*where $\nu = \frac{\sqrt{\sigma^2+2C}}{\sigma}$ and we interpret $I_{-\nu}(z)$ to be $K_\nu(z)$ if $\nu$ is an integer.*

PROOF. Lennox proved in [13] that if the drift $f$ satisfies the given Riccati equation, then the PDE has a Lie algebra of symmetries spanned by

$$\mathbf{v}_1 = \partial_t, \qquad \mathbf{v}_2 = xe^{\sqrt{A}t}\partial_x + \frac{e^{\sqrt{A}t}}{\sqrt{A}}\partial_t - \frac{1}{2\sigma}\left(\sqrt{A}x + f(x) + \frac{B}{\sqrt{A}}\right)e^{\sqrt{A}t}u\partial_u,$$

$$\mathbf{v}_3 = -xe^{-\sqrt{A}t}\partial_x + \frac{e^{-\sqrt{A}t}}{\sqrt{A}}\partial_t - \frac{1}{2\sigma}\left(\sqrt{A}x - f(x) + \frac{B}{\sqrt{A}}\right)e^{-\sqrt{A}t}u\partial_u,$$

$$\mathbf{v}_4 = u\partial_u, \qquad \mathbf{v} = \beta(x,t)\partial_u.$$

Here $\beta$ is an arbitrary solution of the PDE. The symmetries $\mathbf{v}_\beta$ corresponding to adding a solution $\beta$ to the original solution. Plainly these cannot give group invariant solutions. We will solve

(5.7)
$$u_t = \sigma x u_{xx} + f(x)u_x - \mu x u, \qquad x \geq 0,$$
$$u(x,0) = \delta(x-y).$$

We look for symmetries of the form $\mathbf{v} = \sum_{k=1}^4 c_k \mathbf{v}_k$ which preserve the boundary conditions. We require the boundary $x = 0$ to be preserved by the action of $\mathbf{v}$. That is $\mathbf{v}(x) = 0$ when $x = 0$. We also require the boundary $t = 0$ to be preserved, which means that $\mathbf{v}(t) = 0$ when $t = 0$. Finally $u(x,0) = \delta(x-y)$ must be preserved, which requires $\mathbf{v}(u - \delta) = 0$ when $u(x,0) = \delta(x-y)$. These conditions are satisfied if and only if $c_2 = c_3, c_1 = -\frac{2}{\sqrt{A}}c_2$ and $c_4 = 2\frac{(Ay+B)}{\sigma\sqrt{A}}c_2$. This gives the form of the symmetry we need.

We then compute the invariants for the action of

$$\mathbf{v} = c_1\left(-\frac{2}{\sqrt{A}}\mathbf{v}_1 + \mathbf{v}_2 + \mathbf{v}_3 + 2\frac{(Ay+B)}{\sigma\sqrt{A}}\mathbf{v}_4\right).$$

Invariants of the action are given by solving the equation $\mathbf{v}(\eta) = 0$. They are readily found by the method of characteristics and we may take the invariants to be $\eta$ and $v$ where

(5.8) $$\eta = \frac{x}{4\sinh^2(\sqrt{A}t/2)},$$

(5.9) $$u = \exp\left(-\frac{(Bt + F(x) - F(y))}{2\sigma} - \frac{\sqrt{A}(x+y)}{2\sigma\tanh(\sqrt{A}t/2)}\right)$$
$$\times v\left(\frac{x}{4\sinh^2(\sqrt{A}t/2)}\right).$$

Now $v(\eta)$ must satisfy $2\sigma^2\eta^2 v''(\eta) - (C + 2Ay\eta)v(\eta) = 0$ if $u$ is a solution of the PDE. Hence

$$v(\eta) = \sqrt{y\eta}\left(C_1(y)I_\nu\left(\frac{2\sqrt{Ay\eta}}{\sigma}\right) + C_2(y)I_{-\nu}\left(\frac{2\sqrt{Ay\eta}}{\sigma}\right)\right).$$

Since the solution of our initial value problem is a group invariant solution for the action generated by $\mathbf{v}$ the result follows from substitution of $v$ and $\eta$ into the expression for $u$. □

NOTE. The symmetries used in the proof are actually symmetries for any PDE $u_t = \sigma x u_{xx} + f(x)u_x - g(x)u$ when $\sigma x f' - \sigma f + \frac{1}{2}f^2 + 2\sigma x g(x) = \frac{1}{2}Ax^2 + Bx + C$. So we may use the same technique to construct the form of the fundamental solutions for any PDE of this type. We will present the case when $g(x) = \frac{\mu}{x}$ below.

5.1. *The functions $C_1(y), C_2(y)$ and transition densities.* The functions $C_1, C_2$ will depend upon the initial and boundary conditions we may wish to impose on the PDE. The Bessel function $I_{-\nu}(y)$ is not integrable near zero for $\nu \geq 1$, so if we require a fundamental solution which defines solutions of the Cauchy problem for constant initial data, then we will have $C_2 = 0$. It is also this choice with $C_1 = 1$ which is needed if we are to recover the transition density as we let $\mu \to 0$. In the general case, $C_1, C_2$ can be found from the condition that $u(x,t) = \int_0^\infty \phi(y)p(t,x,y)\,dy$ solves the PDE with $u(x,0) = \phi(x)$. For most cases, one may simply take $\phi(x) = 1, \phi(x) = x$ for the initial data, from which we may identify $C_1$ and $C_2$.

Recovering a transition density needs care however. Though our expression for the fundamental solutions are very general, it is possible to have fundamental solutions for the types of PDEs we have been considering which involve terms with delta functions. We have seen such a situation in Example 5.1.

In all save one of the examples we present below, we do indeed obtain the transition density by letting $\mu \to 0$ in the fundamental solution we obtain from Theorem 5.4. However in general, letting $\mu = 0$ will produce a fundamental solution of $u_t = \sigma x u_{xx} + f(x) u_x$, but there is no guarantee that this fundamental solution will be the transition density for the diffusion. We need to be careful to incorporate the behavior at the boundary $y = 0$.

One can usually produce a density from the fundamental solutions provided by Theorem 5.4. Suppose that $p_\mu(t,x,y)$ is the fundamental solution we obtain from an application of Theorem 5.4. We then have $\int_0^\infty \phi(y) \times p_\mu(t,x,y)\,dy \to \phi(x)$ as $t \to 0$. Now introduce the function $h_\mu(t,x) = \int_0^\infty p_\mu(t,x,y)\,dy$. If $\lim_{\mu \to 0} h_\mu(t,x) = 1$, then $p_\mu$ will reduce to a density at $\mu = 0$, since an application of Lebesgue's dominated convergence theorem gives $\int_0^\infty \lim_{\mu \to 0} p_\mu(t,x,y)\,dy = 1$. If $\lim_{\mu \to 0} h_\mu(t,x) \neq 1$ then to produce a density, we let $\overline{p}(t,x,y) = \lim_{\mu \to 0} p_\mu(t,x,y)$ and $h(t,x,y) = \lim_{\mu \to 0} h_\mu(t,x)$. Then $\overline{p}(t,x,y) + (1 - h(t,x,y))\delta(y)$ is a density.

In this case to calculate the expectations we want, we will have to incorporate additional terms involving generalized functions to produce a fundamental solution that reduces to the necessary density. Consider the symmetry solution (5.1), where we choose a stationary solution not left invariant by the symmetry. We make the following observation: $\lim_{t \to 0} U_1(x,t) = 0$. This follows from the fact that the stationary solutions of the PDE (5.7) are of the general form

$$
(5.10) \quad u_0(x) = x^{\beta/2} e^{-(F(x) - \sqrt{A}x)/(2\sigma)} \\
\times \left( c_1 {}_1F_1\left(\alpha, \beta, \frac{\sqrt{A}x}{\sigma}\right) + c_2 \Psi\left(\alpha, \beta, \frac{\sqrt{A}x}{\sigma}\right) \right)
$$

with $\alpha = \frac{1}{2\sigma}(\frac{B}{\sqrt{A}} + \sigma(1 + \sqrt{1 + 2C/\sigma^2}))$, $\beta = 1 + \sqrt{1 + 2C/\sigma^2}$. The estimates ${}_1F_1(a,b,z) = \frac{\Gamma(b)}{\Gamma(a)} e^z z^{a-b}(1 + \mathrm{O}(|z|^{-1}))$ and $\Psi(a,b,z) = z^{-a}(1 + \mathrm{O}(|z|^{-1}))$ for large $|z|$ with $z > 0$ give the result. (See page 504 of [1].) Thus if $q_\mu = p_\mu(t,x,y) + U_1(x,t)\delta(y)$, then

$$
\int_0^\infty \phi(y) q_\mu(t,x,y)\,dy = \int_0^\infty \phi(y) p_\mu(t,x,y)\,dy + \phi(0) U_1(x,t) \to \phi(x)
$$

as $t \to 0$. Hence $q_\mu$ is also a fundamental solution.

We apply this to the diffusion in Example 5.1. An application of Theorem 5.4 leads to the fundamental solution

$$
p_\mu(t,x,y) = \frac{2+ay}{2+ax} \frac{e^{\sqrt{\mu}y}}{e^{\sqrt{\mu}x}} \exp\left\{ \frac{-2\sqrt{\mu}(x + y e^{2\sqrt{\mu}t})}{e^{2\sqrt{\mu}t} - 1} \right\} \\
\times \sqrt{\frac{\mu x}{y}} \frac{I_1(2\sqrt{\mu x y}/\sinh(\sqrt{\mu}t))}{\sinh(\sqrt{\mu}t)}.
$$

If we add $U_1(x,t)\delta(y)$ to this we obtain the necessary fundamental solution found in Example 5.1. We do not have a proof that this method works in all cases, but similar procedures will always allow us to construct the necessary fundamental solutions so that we can compute the desired expectations.

COROLLARY 5.5. *For $n \geq 2$ a fundamental solution of the PDE*

$$u_t = 2xu_{xx} + nu_x - \frac{b^2}{2}u, \qquad x \geq 0,$$

*is*

$$p(t,x,y) = \frac{b}{2\sinh(bt)}\left(\frac{y}{x}\right)^{n/4-1/2}$$
$$\times \exp\left(-\frac{b(x+y)}{2\tanh(bt)}\right) I_{(n-2)/2}\left(\frac{b\sqrt{xy}}{\sinh(bt)}\right). \quad (5.11)$$

PROOF. We simply need to check that $f$ satisfies the Riccati equation and determine $A, B, C$. We choose $C_1 = 1, C_2 = 0$. We further note that as $b \to 0$, the fundamental solution converges to the transition density for an $n$ dimensional squared Bessel process. One may find other fundamental solutions which do not reduce to the transition density. □

We now have a very easy proof of a well known result for squared Bessel processes, which is usually proved using martingale methods.

COROLLARY 5.6. *Let $X = \{X_t : t \geq 0\}$ be a squared Bessel processes, where*

$$dX_t = n\, dt + 2\sqrt{X_t}\, dW_t, \qquad n \geq 2.$$

*Then*

$$E_x\left(e^{-\lambda X_t - (b^2/2)\int_0^t X_s\, ds}\right)$$
$$= \frac{\exp(-(xb/2)(1 + 2\lambda b^{-1}\coth(bt))/(\coth(bt) + 2\lambda b^{-1}))}{(\cosh(bt) + 2\lambda b^{-1}\sinh(bt))^{n/2}}. \quad (5.12)$$

PROOF. The expectation is given by

$$E_x\left(e^{-\lambda X_t - (b^2/2)\int_0^t X_s\, ds}\right) = \int_0^\infty e^{-\lambda y} p(t,x,y)\, dy, \quad (5.13)$$

where

$$p(t,x,y) = \frac{b}{2\sinh(bt)}\left(\frac{y}{x}\right)^{n/4-1/2}$$
$$\times \exp\left(-\frac{b(x+y)}{2\tanh(bt)}\right) I_{(n-2)/2}\left(\frac{b\sqrt{xy}}{\sinh(bt)}\right). \quad (5.14)$$

Evaluation of the integral is routine. See the tables in section 6.6 of [11]. □

5.2. *The zero coupon bond price in the CIR model.* Using our results, it is easy to compute the well-known zero coupon bond price in the CIR model. We have $A = b^2 + 4\mu\sigma, B = -ab, C = \frac{1}{2}a^2 - a\sigma$. With $C_1(y) = 1, C_2(y) = 0$ and $F(x) = a \ln x - bx$ we recover the necessary fundamental solution. We leave this to the reader. Let us now consider an example where we have to include a delta function term.

EXAMPLE 5.2. We are interested in the process $X$ where
$$dX_t = 2X_t \tanh(X_t)\, dt + \sqrt{2X_t}\, dW_t, \qquad X_0 = x.$$

Using the symmetry (5.1), one may check that the transition density of this process is

(5.15) $q(t,x,y) = \dfrac{1}{\sinh t}\dfrac{\cosh y}{\cosh x}\exp\left\{-\dfrac{x+y}{\tanh t}\right\}\left(\sqrt{\dfrac{x}{y}}I_1\left(\dfrac{2\sqrt{xy}}{\sinh t}\right) + \sinh t\, \delta(y)\right).$

Using our results we find a fundamental solution of

(5.16) $$u_t = xu_{xx} + 2x\tanh(x)u_x - \mu x u$$

is

(5.17) $p(x,y,t) = \dfrac{1}{\sinh kt}\dfrac{\cosh y}{\cosh x}\exp\left\{-\dfrac{k(x+y)}{\tanh kt}\right\}\sqrt{\dfrac{kx}{y}}I_1\left(\dfrac{2\sqrt{kxy}}{\sinh kt}\right).$

Here $k = \sqrt{1+\mu}$. This does not reduce to the transition density as $\mu \to 0$. However, $u_0(x) = e^{-\sqrt{1+\mu}x}\operatorname{sech} x$ is a stationary solution of this PDE. If we apply the symmetry (5.1) to this solution, we find a new solution $U_\varepsilon(x,t)$ which has the property that $U_1(x,t) = \frac{1}{\cosh x}\exp\{-\frac{kx}{\tanh kt}\}$.

Now $q(x,y,t) = p(x,y,t) + \delta(y)U_1(x,t)$ is also a fundamental solution and this does reduce to the transition density. From which we may compute

$$\mathbb{E}_x(e^{-\lambda X_t - \mu \int_0^t X_s\, ds})$$
$$= \int_0^\infty e^{-\lambda y}q(x,y,t)\, dy$$
$$= U_1(x,t) + \dfrac{e^{-x/\tanh(kt)}}{2}(e^{kx\operatorname{Csch}(kt)/(\cosh(kt)+(-1+\lambda)\sinh(kt))}$$
$$+ e^{kx\operatorname{Csch}(kt)/(\cosh(kt)+(1+\lambda)\sinh(kt))} - 2).$$

5.3. *Calculating expectations of the form $E_x(e^{-\lambda X_t - \mu \int_0^t \frac{ds}{X_s}})$.* In Section 3 we calculated the expectation $\xi_{\lambda,\mu/x}(X_t)$ for the case when $\gamma = 1$ and $\sigma x f' - \sigma f + \frac{1}{2}f^2 = Ax + B, A \geq 0$. We now show how to compute these expectations when $f$ satisfies (2.11).

THEOREM 5.7. *Suppose that $\sigma x f' - \sigma f + \frac{1}{2}f^2 = \frac{1}{2}Ax^2 + Bx + C$, $A > 0$. Then the PDE*

$$u_t = \sigma x u_{xx} + f(x)u_x - \frac{\mu}{x}u, \qquad \mu \geq 0, \tag{5.18}$$

*has a fundamental solution of the form*

$$
\begin{aligned}
&p(t,x,y) \\
&= \frac{\sqrt{A}e^{(F(y)-F(x))/(2\sigma)}}{2\sigma \sinh(\sqrt{A}t/2)}\sqrt{\frac{x}{y}}\exp\left\{-\frac{Bt}{2\sigma} - \frac{\sqrt{A}(x+y)}{2\sigma\tanh(\sqrt{A}t/2)}\right\} \\
&\quad \times \left(C_1(y)I_\nu\left(\frac{\sqrt{Axy}}{\sigma\sinh(\sqrt{A}t/2)}\right) + C_2(y)I_{-\nu}\left(\frac{\sqrt{Axy}}{\sigma\sinh(\sqrt{A}t/2)}\right)\right),
\end{aligned}
\tag{5.19}
$$

*in which $F'(x) = f(x)/x$ and $\nu = \frac{\sqrt{2C+4\mu\sigma+\sigma^2}}{\sigma}$ and we interpret $I_{-\nu}(z)$ to be $K_\nu(z)$ if $\nu$ is an integer.*

PROOF. The proof is similar to the previous result. The invariants are the same, but the PDE in this case reduces to the ODE

$$4\sigma^2\eta^2 v''(\eta) - (2C + 4Ay\eta + 4\mu\sigma)v(\eta) = 0. \qquad \square$$

The following result for Cox–Ingersoll–Ross process is an easy corollary.

COROLLARY 5.8. *The PDE $u_t = \sigma x u_{xx} + (a - bx)u_x - \frac{\mu}{x}u$ with $\mu \geq 0$, $a, b > 0$ has a fundamental solution*

$$
\begin{aligned}
p(t,x,y) &= \frac{b}{2\sigma\sinh(bt/2)}\left(\frac{y}{x}\right)^{a/(2\sigma)-1/2} \\
&\quad \times \exp\left(\frac{b}{2\sigma}\left(at + (x-y) - \frac{x+y}{\tanh(bt/2)}\right)\right) \\
&\quad \times I_\nu\left(\frac{b\sqrt{xy}}{\sigma\sinh(bt/2)}\right).
\end{aligned}
\tag{5.20}
$$

*Here $\nu = \frac{1}{\sigma}\sqrt{(a-\sigma)^2 + 4\mu\sigma}$.*

COROLLARY 5.9. *Let*

$$k = \frac{a}{2\sigma}, \qquad \alpha = \frac{b}{2\sigma}\left(1 + \coth\left(\frac{bt}{2}\right)\right) + \lambda, \qquad \beta = \frac{b\sqrt{x}}{2\sigma\sinh(bt/2)}.$$

*and $M_{s,r}(z)$ be the Whittaker functions of the first kind. For the CIR process $dX_t = (a - bX_t)\,dt + \sqrt{2\sigma X_t}\,dW_t$ we have*

$$E_x(e^{-\lambda X_t - \mu \int_0^t \frac{ds}{X_s}}) = \frac{\Gamma(k+\nu/2+1/2)}{\Gamma(\nu+1)}\beta x^{-k}$$

$$\text{(5.21)} \quad \times \exp\left(\frac{b}{2\sigma}\left(at + x - \frac{x}{\tanh(bt/2)}\right)\right)$$
$$\times \frac{1}{\beta\alpha^k} e^{\beta^2/(2\alpha)} M_{-k,\nu/2}\left(\frac{\beta^2}{\alpha}\right).$$

PROOF. Observe that the fundamental solution in Corollary 5.8 reduces to the transition density of a CIR process as $\mu \to 0$. So we have

$$\text{(5.22)} \quad E_x(e^{-\lambda X_t - \mu \int_0^t \frac{ds}{X_s}}) = \int_0^\infty e^{-\lambda y} p(t,x,y)\, dy.$$

The result follows from the fact that

$$\int_0^\infty y^{k-1/2} e^{-\alpha y} I_{2\gamma}(2\beta\sqrt{y})\, dy = \frac{\gamma(k+\gamma+1/2)}{\Gamma(2\gamma+1)} \frac{1}{\beta\alpha^k} e^{\beta^2/(2\alpha)} M_{-k,\gamma}\left(\frac{\beta^2}{\alpha}\right),$$

which is formula 6.643.2 of [11]. $\square$

5.4. *Laplace transforms of joint densities for* $(X_t, \int_0^t X_s\, ds, \int_0^t \frac{ds}{X_s})$. By similar means to the above, we may find fundamental solutions of the PDE $u_t = \sigma x u_{xx} + f(x)u_x - (\frac{\nu}{x} + \mu x)u$.

THEOREM 5.10. *Suppose that $f$ is a solution of the Riccati equation $\sigma x f' - \sigma f + \frac{1}{2}f^2 + 2\sigma\nu + 2\sigma\mu x^2 = \frac{1}{2}Ax^2 + Bx + C, A > 0$. Then the PDE*

$$\text{(5.23)} \quad u_t = \sigma x u_{xx} + f(x)u_x - \left(\frac{\nu}{x} + \mu x\right)u, \qquad \mu > 0, \nu > 0,$$

*has a fundamental solution of the form*

$$p(t,x,y) = \frac{\sqrt{Axy}}{2\sigma \sinh(\sqrt{A}t/2)} e^{-(Bt + \sqrt{A}(x+y)\coth((\sqrt{A}t/2)) + F(x) - F(y))/(2\sigma)}$$

$$\text{(5.24)} \quad \times \left(C_1(y) I_{\sqrt{\sigma^2 + 2C}/\sigma}\left(\frac{\sqrt{Axy}}{\sigma \sinh(\sqrt{A}t/2)}\right)\right.$$
$$\left. + C_2(y) I_{-\sqrt{\sigma^2 + 2C}/\sigma}\left(\frac{\sqrt{Axy}}{\sigma \sinh(\sqrt{A}t/2)}\right)\right).$$

*As usual $F'(x) = f(x)/x$ and $I_{-d}(z) = K_d(z)$ if $d$ is an integer.*

PROOF. The proof is similar to the previous cases, with the PDE having the same infinitesimal symmetries and invariants. $\square$

Again, we will normally have $C_1 = 1, C_2 = 0$.

EXAMPLE 5.3. We consider a squared Bessel process of dimensions $n$. We apply the previous result with $A = 4b^2$, $B = 0$, $C = 4\nu - 2n + \frac{1}{2}n^2$, $C_1 = 1$, $C_2 = 0$. The PDE $u_t = 2xu_{xx} + nu_x - (\frac{b^2}{2}x + \frac{\nu}{x})u$ has a fundamental solution

$$p(t,x,y) = \frac{b}{2\sinh(bt)} e^{-b(x+y)/(2\tanh(bt))} \left(\frac{y}{x}\right)^{(n-2)/4} I_{\sqrt{(n-2)^2+8\nu}/2}\left(\frac{b\sqrt{xy}}{\sinh(bt)}\right).$$

This reduces to the necessary density and we find

$$E_x\left(e^{-\lambda X_t - (b^2/2)\int_0^t X_s\,ds - \nu \int_0^t \frac{ds}{X_s}}\right)$$

$$= \int_0^\infty e^{-\lambda y} p(t,x,y)\,dy$$

$$= e^{-bx/(2\tanh bt)} \frac{\Gamma(\alpha)}{\Gamma(\beta)} \frac{b^{a/2}(xe^{bt})^\gamma (e^{2bt}-1)^{-\gamma}}{(\cosh bt + (2\lambda/b)\sinh bt)^\delta}$$

$$\times {}_1F_1\left(\alpha, \beta, \frac{b^2 x\,\text{csch}(bt)}{2b\cosh(bt) + 4\lambda\sinh(bt)}\right),$$

where $a = \sqrt{(n-2)^2 + 8\nu}$, $\delta = \frac{1}{4}(2+a+n)$, $\gamma = \frac{1}{4}(2+a-n)$, $\alpha = \frac{1}{4}(a+n+2)$ and $\beta = \frac{a+2}{2}$.

**5.5. An example in the $\gamma = 0$ case.** We prove the following result.

THEOREM 5.11. *Let $f$ be a solution of the Riccati equation $\sigma f' + \frac{1}{2}f^2 = \frac{C}{x^2} + \frac{B}{2} + \frac{1}{8}(A - 16\mu\sigma)x^2$, $A > 0$. Then the PDE*

(5.25) $$u_t = \sigma u_{xx} + f(x)u_x - \mu x^2 u, \qquad x > 0,$$

*has a fundamental solution of the form*

(5.26)
$$p(t,x,y) = \frac{\sqrt{Axy}}{2\sinh(\sqrt{A}t/2)}$$
$$\times \exp\left(-\frac{Bt}{4\sigma} + \frac{F(y)-F(x)}{2\sigma} - \frac{\sqrt{A}(x^2+y^2)}{8\sigma\tanh(\sqrt{A}t/2)}\right)$$
$$\times \left(C_1(y)I_\beta\left(\frac{\sqrt{A}xy}{4\sigma\sinh(\sqrt{A}t/2)}\right)\right.$$
$$\left. + C_1(y)I_{-\beta}\left(\frac{\sqrt{A}xy}{4\sigma\sinh(\sqrt{A}t/2)}\right)\right).$$

*Here $F' = f$, $\beta = \frac{1}{2}\sqrt{1 + 2C/\sigma^2}$ and we interpret $I_{-\beta}$ to be $K_\beta$ if $\beta$ is an integer.*

PROOF. Applying Lie's algorithm shows that if $f$ satisfies the given Riccati equation, then the finite dimensional part of the Lie algebra of symmetries is spanned by $\mathbf{v}_1 = \partial_t, \mathbf{v}_4 = u\partial_u$,

$$\mathbf{v}_2 = \frac{\sqrt{A}}{2}xe^{\sqrt{A}t}\partial_x + e^{\sqrt{A}t}\partial_t - \left(\frac{A}{8\sigma}x^2 + \frac{\sqrt{A}}{4\sigma}xf(x) + m\right)e^{\sqrt{A}t}u\partial_u,$$

$$\mathbf{v}_3 = -\frac{\sqrt{A}}{2}xe^{-\sqrt{A}t}\partial_x + e^{-\sqrt{A}t}\partial_t - \left(\frac{A}{8\sigma}x^2 - \frac{\sqrt{A}}{4\sigma}xf(x) - n\right)e^{-\sqrt{A}t}u\partial_u,$$

$m = \frac{\sqrt{A}+B/\sigma}{4}$, $n = \frac{\sqrt{A}-B/\sigma}{4}$. Proceeding as previously, we find that the fundamental solution has the form

$$p(t,x,y) = e^{-2Bt+\sqrt{A}(x^2+y^2)\coth(\sqrt{A}t/2)+4\sigma \ln(\sinh(\sqrt{A}t/2))/(8\sigma)} v\left(\frac{x}{\sinh(\sqrt{A}t/2)}\right).$$

Where $16\eta^2\sigma^2 v''(\eta) - (Ay^2\eta^2 + 8C)v(\eta) = 0$. The result follows. □

As previously, we will normally take $C_1 = 1$, $C_2 = 0$.

EXAMPLE 5.4. We consider a radial Ornstein–Uhlenbeck process. The SDE is $dX_t = (\frac{a}{X_t} + bX_t) + \sqrt{2}\,dW_t$, $X_0 = x$. For simplicity we take $a > 1/2$. The fundamental solution is found by the previous theorem to be

$$
\begin{aligned}
p(t,x,y) &= y\left(\frac{y}{x}\right)^{\nu-1}\frac{\alpha}{\sinh(\alpha t)} \\
&\quad \times e^{-\nu t - \alpha(x^2+y^2)/(4\tanh(\alpha t)) - (b/4)(x^2-y^2)} I_{\nu-1}\left(\frac{\alpha xy}{2\sinh(\alpha t)}\right),
\end{aligned}
$$
(5.27)

$\alpha = \sqrt{b^2 + 4\mu}$ and $\nu = \frac{1}{2}(a+1)$. Note that as $\mu \to 0$, this reduces to the transition density for a radial Ornstein–Uhlenbeck process. (See [9] and [10].) From which we find for example that

$$E_x(e^{-\lambda X_t^2 - \mu \int_0^t X_s^2\,ds}) = \int_0^\infty e^{-\lambda y^2} p(t,x,y)\,dy$$

$$= 2e^{-bx^2/4 + \alpha(\alpha - (b-4\lambda)\coth(t\alpha))x^2/(4(b-4\lambda - \alpha\coth(t\alpha))) - \nu t}$$

$$\times \left(\cosh(t\alpha) - \frac{(b-4\lambda)\sinh(t\alpha)}{\alpha}\right)^{-\nu}.$$

**6. Computable functionals for a given drift.** Observe that the roles of $g(x)$ and $f(x)$ can be swapped in the sense that we can allow $f(x)$ to be a fixed drift and we can then determine the functions $g$ for which the PDEs (1.3) has a nontrivial Lie group of symmetries. Thus, even though we cannot determine the transition density for all drifts $f$, we can always compute specific functionals $\xi_{\lambda,g}(X_t)$ for a given drift.

EXAMPLE 6.1. A simple example illustrates the essential idea. Consider the drift $f(x) = a - b\sqrt{x}$. Suppose $\sigma = \gamma = 1$. We take $A, B$ to be arbitrary, but for simplicity we assume that $A > 0, B > 0$. If $xf' - f + 1/2f^2 + 2xg(x) = Ax + B$ then clearly the function $g$ must be of the form

$$g(x) = \frac{A - b^2/2}{2} + \frac{a - a^2/2 + B}{2x} + \frac{ab - b/2}{2\sqrt{x}}.$$

From the stationary solution $u_0(x) = x^{(1-a)/2} e^{bx/2} I_{\sqrt{1+2B}}(\sqrt{2Ax})$ we find that

$$(6.1) \quad U_\lambda(x,t) = \frac{x^{1/2-a/2}}{1+\lambda t} \exp\left\{\frac{1}{2}bx - \frac{\lambda(x + At^2/2)}{1+\lambda t}\right\} I_{\sqrt{1+2B}}\left(\frac{\sqrt{2Ax}}{(1+\lambda t)}\right)$$

is a solution of the PDE $u_t = xu_{xx} + (a - b\sqrt{x})u_x - g(x)u$. Inverting the Laplace transform to obtain the fundamental solution gives

$$(6.2) \quad p(t,x,y) = \frac{1}{t}\left(\frac{x}{y}\right)^{(1-a)/2} \exp\left\{\frac{b(x-y) - At}{2} - \frac{x+y}{t}\right\} \times I_{\sqrt{1+2B}}\left(\frac{2\sqrt{xy}}{t}\right).$$

One can check that this is the correct fundamental solution and so that if the diffusion $X_t$ satisfies the SDE $dX_t = (a - b\sqrt{X_t})\,dt + \sqrt{2X_t}\,dW_t$, then

$$E_x(e^{-\lambda X_t + \int_0^t g(X_s)\,ds})$$
$$= 2^a (2xt)^{(\mu+1-a)/2} e^{(-At + (b-2/t)x)/2} (2 + t(b+2\lambda))^{(-1-a-\mu)/2}$$
$$\times \frac{\Gamma((1+a+\mu)/2)}{\Gamma(1+\mu)} {}_1F_1\left(\frac{1+a+\mu}{2}, 1+\mu, \frac{2x}{t(2+t(b+2\lambda))}\right),$$

where $\mu = \sqrt{1+2B}$. This is not enough to deduce the transition density except in the case $a = \frac{1}{2}$. However these calculations can be carried out for any one dimensional diffusion and can be potentially quite useful.

**Acknowledgments.** The authors wish to thank Eckhard Platen, the Australian Research Council and St. George's Bank Ltd. for their support during the completion of this work.

Department of Mathematical Sciences
University of Technology, Sydney
P.O. Box 123, Broadway
New South Wales 2007
Australia
E-mail: Mark.Craddock@uts.edu.au
        Kelly.A.Lennox@student.uts.edu.au